\documentclass[a4paper,12pt,draft]{amsart}
\usepackage{fullpage}
\usepackage{amsmath}
\usepackage{amsfonts}
\usepackage{amssymb}
 \def\o{\omega}

\def\e{\epsilon}
\def\w{\wedge}

\def\va{\varphi}

\def\c{\centerline}
\begin{document}
\title{\bf Convergence in capacity on compact K\"ahler manifolds}
\author{S\l awomir Dinew, Ph\d{a}m Ho\`ang Hi\d{\^{e}}p}
\maketitle
\begin{abstract}
 The aim of this  note is to study the convergence in capacity for functions in the class $\mathcal E(X,\o )$. We obtain several stability theorems. Some of these are (optimal) generalizations of results of Xing [Xi2], while others are new.
\end{abstract}
\vskip 0,5 cm
\noindent
\centerline{\bf INTRODUCTION}
\medskip

\noindent
The notion of relative capacity $C_n$ was introduced by Bedford and Taylor (cf. [BT2]). The convergence with respect to $C_n$ in the setting of domains in $\bf C^n$ was investigated by Xing and Cegrell (see [Xi1], [Ce3]). This topic has attracted much interest, since the complex Monge-Amp\`ere operator is continuous with respect to such a convergence (while continuity fails for, say, convergence in $L^p, 1<p<\infty$, see [CK1]). Recently Ko{\l}odziej (see [Ko1]) introduced the capacity $C_{X,\omega}$ on a compact K\"ahler manifold $(X,\omega)$, which is modelled on $C_n$.
 In particular he proved that $C_{X,\omega}$\ is locally equivalent to $C_n$. This capacity was also studied by Guedj and Zeriahi in [GZ1].  In
 [GZ2] they introduced the new Cegrell class $\mathcal E(X, \omega)$ of $\omega$-psh functions for which the complex Monge-Amp\`{e}re operator is
 well-defined. Roughly speaking, it is the largest class of $\omega$-psh functions where the Monge-Amp\`ere operator behaves like it does in the
 bounded functions setting. 

The aim of the present note is to study the convergence in capacity $C_{X,\omega}$ in the class $\mathcal E(X,\omega)$ (Section 2). Thus this note might be seen as a K\"ahler manifold counterpart of [CK2]. Some of our results are known for subclasses of $\mathcal E(X,\omega)$, (especially for bounded or $\mathcal E^1$ functions - see [Ko2] and [Xi2]), or under additional technical assumptions (cf. [Xi2]). By applying some new results borrowed from [Di2] and [Di3] we are able to obtain the sharp results and relax those conditions.

Some of our results, however, (as well as the techniques used) are completely new, thus these might be of independent interest.

{\bf Acknowledgments.} The first named author was partially supported by Polish ministerial grant N N 201 271135.
\vskip 0,5 cm

\c {\bf 1. PRELIMINARIES}
\medskip
\noindent
First we recall some elements of pluripotential theory that will be used throughout the paper. All this can be found in [Bl], [BT1-2], [Ce1-3], [CGZ], [Di1-3], [GZ1-2], [H1-2], [H\"o], [KH], [Ko1] and [Xi1-2].

{\bf 1.1.}  Let $X$ be a compact K\"ahler manifold with a fundamental form $\omega=\omega_X$ with $\int\limits_X\omega^n=1$. An upper semicontinuous function $\varphi :\ X\to [-\infty,+\infty)$ is called $\omega$-plu\-ri\-sub\-har\-monic ($\omega$-psh for short) if $\varphi\in L^1(X)$ and $\omega_\varphi :=\omega + dd^c\varphi\geq 0$, where the inequality is understood in the sense on currents. By $PSH(X,\omega )$ (resp. $PSH^-(X,\omega )$) we denote the set of $\omega$-psh (resp. negative $\omega$-psh) functions on $X$.

{\bf 1.2.} In [Ko1] Ko{\l}odziej introduced the capacity $C_{X,\o }$ on $X$ by
$$C_X(E)=C_{X,\o }(E)=\sup\{\int\limits_E\o_\va^n:\ \va\in\text{PSH}(X,\o),\ -1\leq\va\leq 0\}$$
for any Borel set $E\subset X$, where $\o_\va^n=(\o +dd^c\va )^n$ and $n=dim X$.

We refer to [Ko1], [GZ1] for more informations about this capacity.

{\bf 1.3.} The following class of $\o$-psh fuctions was introduced by Guedj and Zeriahi in [GZ2]:
$$\mathcal E(X,\o )=\{\va\in PSH(X,\o ):\ \lim\limits_{j\to\infty }\int\limits_{\{\va >-j\}}\o_{\max (\va ,-j)}^n=\int\limits_X\o^n=1\}.$$
Intuitively it consists of those functions which have mild singularities (or do not have $-\infty$ poles at all), so that no Monge-Amp\`ere mass concentrates near those poles. In particular any such function must have zero Lelong numbers everywhere.

Let us also define
$$\mathcal E^- (X,\o )=\mathcal E(X,\o )\cap PSH^-(X,\o ).$$

We refer to [GZ2] for all the properties of functions from $\mathcal E(X,\o)$.

{\bf 1.4.} Let $u_j,u\in\text{PSH}(X,\o )$. We say that $\{u_j\}$ converges to $u$ in $C_X$ if
$$C_X(\{|u_j-u|>\delta\})\to 0$$
as $j\to\infty$, for every (fixed) $\delta >0$.

{\bf 1.5.} A family of positive measures $\{\mu_\alpha\}$ on $X$ is said to be uniformly absolutely continuous with respect to $C_X$-capacity if for every $\epsilon >0$ there exists $\delta >0$ such that for each Borel subset $E\subset X$ satisfying $C_X(E)<\delta$ the inequality $\mu_\alpha (E)<\epsilon$ holds for all $\alpha$. We denote this by $\mu_\alpha\ll C_X$ uniformly for $\alpha$. In particular all such measures must vanish on pluripolar sets.

\bigskip
Next we state some well-known results needed for our work. We shall sketch some of the proofs for the sake of completeness:

{\bf 1.6. Proposition.} {\sl Let $u\in\mathcal E(X,\o )$, $v\in\text{PSH}(X,\o )$. Then
$$\int\limits_{\{u\leq v\}}\o_{\max (u,v)}^k\w T=\int\limits_{\{u\leq v\}}\o_u^k\w T,$$
for all $1\leq k\leq n$ and $T=\o_{\varphi_1}\w ... \w\o_{\varphi_{n-k}}$ with $\varphi_1,...,\varphi_{n-k}\in\mathcal E(X,\o )$.}

{\sl Proof.} By Corollary 1.7 in [GZ2] we get
$$\aligned&\int\limits_{\{u\leq v\}}\o_{\max (u,v)}^k\w T=\int\limits_{X}\o_{\max (u,v)}^k\w T-\int\limits_{\{u>v\}}\o_{\max (u,v)}^k\w T\\
&=\int\limits_{X}\o_u^k\w T-\int\limits_{\{u>v\}}\o_u^k\w T=\int\limits_{\{u\leq v\}}\o_u^k\w T.\endaligned$$

{\bf 1.7. Proposition.} {\sl Let $u\in\text{PSH}^-(X,\o )$. Then for $t\geq 0$
$$C_X(\{u<-t\})\leq \frac {|\sup\limits_X u|+c} t,$$
where the positive constant $c$ does not depend on $u$.}

{\sl Proof.} 
$$\aligned &C_X(\{u<-t\})=\sup\{\int\limits_{\{u<-t\}}\o_{\va}^n:\ \va\in PSH(X,\o),\ -1\leq\va\leq 0\}\\
&\leq\sup\{\int\limits_{\{u<-t\}}\frac {|u|} t\o_{\va}^n:\ \va\in PSH(X,\o),\ -1\leq\va\leq 0\}\\&
\leq \frac 1 t\sup\{\int\limits_X |u|\o_{\va}^n:\ \va\in PSH(X,\o),\ -1\leq\va\leq 0\}\\
&\leq \frac 1 t [\int\limits_X |u| \o^n+n]\leq \frac {|\sup\limits_X u|+c} t,\endaligned$$
where we have used Proposition 1.7 and Corollary 2.3 from [GZ1] in the last two inequalities.

{\bf 1.8. Proposition.} {\sl Let $u_j,\ u\in PSH(X,\o )$. Then the following two statements are equivalent:

i) $u_j - u\to 0$ in $C_X$;

ii) $\sup\limits_{j\geq 1}|\sup\limits_X u_j| < + \infty$ and $\max (u_j, -t) - \max (u, -t)\to 0$ in $C_X$ for every $t>0$.}

{\sl Proof.} i) $\Rightarrow$ ii) By a result from [H\"o] we have that $\int\limits_X u_j\o^n\rightarrow\int\limits_X u\o^n>-\infty,\ j\rightarrow\infty$. Thus we obtain
$$\inf\limits_{j\geq 1}\sup\limits_X u_j\geq\inf\limits_{j\geq 1}\int\limits_X u_j\o^n >-\infty .$$
Since $|u_j-u|\geq |\max (u_j, -t)-\max (u, -t)|$ we get
$$C_X (\{|\max (u_j, -t)-\max (u, -t)|>\delta\})\leq C_X (\{|u_j-u|>\delta\})\to 0,$$
for all $\delta >0$, $t>0$.

ii) $\Rightarrow$ i) We can assume that $u_j,u\leq 0$. Set 
$$M=\sup\limits_{j\geq 1}|\sup\limits_X u_j|+|\sup\limits_X u|$$
By Proposition 1.7 we have
$$\aligned &C_X(\{|u_j-u|>3\delta\})\leq C_X(\{|u_j-\max (u_j,-t)|>\delta\})+C_X(\{|\max (u,-t)-u|>\delta\})\\
&\ +C_X(\{|\max (u_j,-t)-\max (u,-t)|>\delta\})\leq C_X(\{u_j<-t\})+C_X(\{u<-t\})\\
&\ +C_X(\{|\max (u_j,-t)-\max (u,-t)|>\delta\})\leq\frac {|\sup\limits_X u_j|+|\sup\limits_X u|+2c} t\\
&+C_X(\{|\max (u_j,-t)-\max (u,-t)|>\delta\})\leq\frac {M+2c} t\\
&+C_X(\{|\max (u_j,-t)-\max (u,-t)|>\delta\})\endaligned$$
for any $t>t_0>0$, $\delta >0$. Hence
$$\varlimsup\limits_{j\to\infty} C_X(\{|u_j-u|>3\delta\})\leq\frac {M+2c} t.$$
Finally by letting $t\to+\infty$ we get
$$\lim\limits_{j\to\infty} C_X(\{|u_j-u|>3\delta\})=0$$
for every $\delta >0$.

{\bf 1.9. Proposition.} {\sl Let $u_j\in\mathcal E^-(X,\o )$. Then the following two statements are equivalent:

i) $\inf\limits_{j\geq 1}\sup\limits_X u_j>-\infty$ and $\o_{u_j}^n\ll C_X$ uniformly for $j\geq 1$;

ii) $\lim\limits_{t\to +\infty}\varlimsup\limits_{j\to\infty} \int\limits_{\{u_j\leq -t\}} \o_{u_j}^n=0$.} 

{\sl Proof.} i) $\Rightarrow$ ii) This is a direct application of Proposition 1.7.

ii) $\Rightarrow$ i) We claim that $\inf\limits_{j\geq 1}\sup\limits_X u_j>-\infty$. Suppose on contrary that $\inf\limits_{j\geq 1}\sup\limits_X u_j=-\infty$. Hence
$$\varlimsup\limits_{j\to\infty} \int\limits_{\{u_j\leq -t\}} \o_{u_j}^n=\varlimsup\limits_{j\to\infty} \int\limits_X \o_{u_j}^n=1,$$
but this contradicts
$$\lim\limits_{t\to +\infty}\varlimsup\limits_{j\to\infty} \int\limits_{\{u_j\leq -t\}} \o_{u_j}^n=0.$$
Let us fix $\e>0$. We choose $t_0$ and then $j_0=j_0(t_0)>1$ such that
$$\int\limits_{\{u_j\leq -t_0\}} \o_{u_j}^n<\e,$$
for all $j\geq j_0$. For each Borel set $E\subset X$ we have
$$\aligned&\int\limits_E \o_{u_j}^n=\int\limits_{E\cap\{u_j\leq -t_0\}} \o_{u_j}^n+\int\limits_{E\cap\{u_j>-t_0\}} \o_{u_j}^n\leq\int\limits_{\{u_j\leq -t_0\}} \o_{u_j}^n+\int\limits_{E\cap\{u_j>-t_0\}} \o_{\max (u_j,-t_0)}^n\\
&\leq\int\limits_{\{u_j\leq -t_0\}} \o_{u_j}^n+t_0^nC_X(E)\leq\e+t_0^nC_X(E)\endaligned$$
for all $j\geq j_0$. On the other hand, since $\sum\limits_{k=1}^{j_0}\o_{u_k}^n\ll C_X$ we can choose $\delta_1>0$ such that $\sum\limits_{k=1}^{j_0}\o_{u_k}^n (E)<\e$ for all Borel sets $E\subset X$ with $C_X(E)<\delta_1$. Hence
$\o_{u_j}^n (E)<2\e$ for all $j\geq 1$ and all  Borel sets $E\subset X$, such that $C_X(E)<\delta=\min (\delta_1,\frac {\e}{t_0^n})$.

\bigskip
The next proposition is a modified version of Lemma 2 in [Xi2]:

{\bf 1.10. Proposition.} {\sl Let $u_j,v_j\in\mathcal E^-(X,\o )$ be such that $u_j\geq v_j$ for $j\geq 1$. Assume that $\o_{v_j}^n\ll C_X$ uniformly for $j\geq 1$. Let also $\inf\limits_{j\geq 1}\sup\limits_X v_j>-\infty$. Then $\o_{u_j}^n\ll C_X$ uniformly for $j\geq 1$.}

{\sl Proof.} By Theorem 1.5 in [GZ2] we have
$$\int\limits_{\{u_j<-2t\}} \o_{u_j}^n\leq 2^n\int\limits_{\{v_j<\frac {u_j} 2-t\}} \o_{\frac {u_j} 2}^n\leq 2^n\int\limits_{\{v_j<\frac {u_j} 2-t\}} \o_{v_j}^n\leq 2^n\int\limits_{\{v_j<-t\}} \o_{v_j}^n,$$
for every $t>0$. By Proposition 1.9 we obtain $\o_{u_j}^n\ll C_X$ uniformly for $j\geq 1$.

We remark that the proof shows that the result still folds for any family of functions $u_{\alpha}$\ and $v_{\alpha},\ \alpha\in\Lambda$ instead of just sequences.

{\bf 1.11. Proposition.} {\sl Let $u_j\in\mathcal E(X,\o )$, $u\in\text{PSH} (X,\o )$ be such that $u_j\to u$ in $C_X$. Then the following two statements are equivalent:

i) $u\in\mathcal E (X,\o )$;

ii) $\o_{u_j}^n\ll C_X$ uniformly for $j\geq 1$.}

{\sl Proof.} i) $\Rightarrow$ ii) By Proposition 1.6 we have
$$\aligned&\int\limits_{\{u_j\leq -t\}} \o_{u_j}^n=\int\limits_{\{u_j\leq -t\}} \o_{\max (u_j,-t)}^n\leq\int\limits_{\{u\leq -t+1\}} \o_{\max (u_j,-t)}^n+\int\limits_{\{|u_j-u|>1\}} \o_{\max (u_j,-t)}^n\\
&\leq\int\limits_{\{u\leq -t+1\}} \o_{\max (u_j,-t)}^n+t^n(\{|u_j-u|>1\}).\endaligned$$
By coupling the quasi-continuity of $u$ (Theorem 3.5 in [BT2]) and the facts that $\o_{\max (u_j,-t)}^n\newline\ll C_X$ uniformly for $j\geq 1$ and $\o_{\max (u,-t)}^n\ll C_X$  uniformly for $t>0$ (see Proposition 1.10), we get 
$$\varlimsup\limits_{t\to +\infty}\varlimsup\limits_{j\to\infty}\int\limits_{\{u_j\leq -t\}} \o_{u_j}^n\leq\varlimsup\limits_{t\to +\infty}\int\limits_{\{u\leq -t+1\}} \o_{\max (u,-t)}^n=0.$$
By Proposition 1.9 we obtain $\o_{u_j}^n\ll C_X$ uniformly for $j\geq 1$.

ii) $\Rightarrow$ i) This is a part of Theorem 3 in [Xi2].

\bigskip
The last proposition is a somewhat stronger version of implication ii) $\Rightarrow$ i):

{\bf 1.12. Proposition.} {\sl Let $u_j\in\mathcal E(X,\o )$, $u\in\text{PSH} (X,\o )$ be such that $u_j\to u$ in $L^1(X)$. Assume that $\o_{u_j}^n\ll C_X$ uniformly for $j\geq 1$. Then $u\in\mathcal E (X,\o )$.}

{\sl Proof.} Set $v_j=\max (u_j, u)\in\mathcal E(X,\o )$. By Hartogs' lemma (see [H\"o])  we have $v_j\to u$ in $C_X$. Thus, by Proposition 1.10, we get that $\o_{v_j}^n\ll C_X$ uniformly for $j\geq 1$. Now, using Proposition 1.11, we obtain $u\in\mathcal E (X,\o )$.
\vskip 0,5 cm
\noindent
\c {\bf 2. CONVERGENCE IN CAPACITY}
Below we prove various stability results regarding convergence in capacity in the class $\mathcal E(X,\omega)$. We begin with a technical result which will be used later on. To authors' knowledge such a theorem is new although some similar ideas can be found in Lemma 2.3 in [H1] and in Theorem 2 in [Xi2].

{\bf 2.1. Theorem.} {\sl Let $u_j,v_j\in\mathcal E(X,\o )$ be such that

i) $\o_{u_j}^n\ll C_X$ uniformly for $j\geq 1$;

ii) $\lim\limits_{j\to \infty} \int\limits_{\{u_j<v_j-\delta\}}\o_{u_j}^n=0$ for all $\delta >0$.

Then $\lim\limits_{j\to\infty}C_X (\{u_j<v_j-\delta\})=0$ for all $\delta >0$.}

{\sl Proof.} We can assume, by adding constants to both $u_j,\ v_j$ if necessary, that $\sup\limits_X u_j =0$ for all $j\geq 1$. Set $u_{jt}:=\max (u_j, -t)$. For each $k=0,...,n$ we will prove inductively that
\begin{equation}\label{1}\lim\limits_{j\to\infty}\sup\{\int\limits_{\{u_j<v_j-\delta\}}\o_{u_j}^{n-k}\w \o_{\va}^k:\ \va\in\text{PSH}(X,\o ),\ -1\leq\va\leq 0\}=0
 \end{equation}

for every $\delta >0$. If $k=0$ then (\ref{1}) holds by assumption.

 Assume that (\ref{1}) holds for $k-1$. We will prove that
$$\lim\limits_{j\to\infty}\sup\{\int\limits_{\{u_j<v_j-3\delta\}}\o_{u_j}^{n-k}\w\o_{\va}^k:\ \va\in PSH(X,\o ),\ -1\leq\va\leq 0\}=0$$
for any $\delta >0$. We fix $t\geq 1$ and $\va\in PSH X,\o )$ satisfying $-1\leq\va\leq 0$. Then, by Theorem 2.3 in [Di3], Corollary 1.7 in [GZ2] and Proposition 1.6 we have
\medskip
\noindent
$$\aligned &\int\limits_{\{u_j<v_j-3\delta\}}\o_{u_j}^{n-k}\w\o_{\va}^k\leq\int\limits_{\{u_j+\frac {\delta} t u_{jt}<v_j+\frac {\delta}{t}\va -2\delta\}}\o_{u_j}^{n-k}\w\o_{\va}^k\\
&\leq\frac {t+\delta} {\delta} \int\limits_{\{u_j+\frac {\delta} t u_{jt}<v_j+\frac {\delta}{t}\va -2\delta\}}\o_{\frac{v_j+\frac {\delta}{t}\va -2\delta}{1+\frac {\delta} t}}\w\o_{u_j}^{n-k}\w\o_{\va}^{k-1}\\
&\leq\frac {t+\delta} {\delta}\int\limits_{\{u_j+\frac {\delta} t u_{jt}< v_j+\frac {\delta}{t}\va -2\delta\}}\o_{\frac {u_j+\frac {\delta} t u_{jt}}{1+\frac {\delta} t}}\w\o_{u_j}^{n-k}\w\o_{\va}^{k-1}\leq\frac {t} {\delta}\int\limits_{\{u_j < v_j-\delta\}}\o_{u_j}^{n-k+1}\w\o_{\va}^{k-1}\\
&+\int\limits_{\{u_j < v_j-\delta\}}\o_{u_{jt}}\w\o_{u_j}^{n-k}\w\o_{\va}^{k-1}\leq\frac {t+\delta } {\delta} \int\limits_{\{u_j < v_j-\delta\}}\o_{u_j}^{n-k+1}\w\o_{\va}^{k-1}\\
&+\int\limits_{\{u_j\leq -t\}}\o_{u_{jt}}\w\o_{u_j}^{n-k}\w\o_{\va}^{k-1}\leq\frac {t+\delta } {\delta} \int\limits_{\{u_j < v_j-\delta\}}\o_{u_j}^{n-k+1}\w\o_{\va}^{k-1}+\int\limits_{\{u_j\leq -t\}}\o_{u_j}^{n-k+1}\w\o_{\va}^{k-1}.\endaligned$$
By the induction hypothesis we get
$$\varlimsup\limits_{j\to\infty}\sup\{\int\limits_{\{u_j<v_j-3\delta\}}\o_{u_j}^{n-k}\w\o_{\va}^k:\ \va\in PSH(X,\o ),\ -1\leq\va\leq 0\}$$
$$\leq\sup\{\int\limits_{\{u_j\leq -t\}}\o_{u_j}^{n-k+1}\w\o_{\va}^{k-1}:\ \va\in PSH (X,\o ),\ -1\leq\va\leq 0,\ j\geq 1\},$$
for every $t\geq 1$. 

Note that $\frac{u_j+\va}2\geq u_j-1$, thus by Proposition 1.10 we get that $\o_{\frac{u_j+\va}2}^n\ll C_X$ uniformly for $j\geq 1$. Expanding the left hand side we obtain 
$$\o_{u_j}^{n-k+1}\w\o_{\va}^{k-1}\ll C_X\ {\rm uniformly\ for}\ j\geq 1.$$
Now by letting $t\to +\infty$ (and using $\o_{u_j}^{n-k+1}\w\o_{\va}^{k-1}\ll C_X$ uniformly for $j\geq 1$),  from Proposition 1.7 we obtain
$$\lim\limits_{j\to\infty}\sup\{\int\limits_{\{u_j<v_j-3\delta\}}\o_{u_j}^{n-k}\w\o_{\va}^k:\ \va\in PSH(X,\o ),\ -1\leq\va\leq 0\}=0.$$

 Theorem 2.1 has various consequences regarding stability of the Monge-Amp\`ere operator. Below we list some of them.

As a direct corollary we obtain the following theorem:

{\bf 2.2. Theorem.} {\sl Let $u_j,v_j\in\mathcal E(X,\o )$ be such that

i) $\o_{u_j}^n+\o_{v_j}^n\ll C_X$ uniformly for $j\geq 1$;

ii) $\lim\limits_{j\to \infty}[ \int\limits_{\{u_j<v_j-\delta\}}\o_{u_j}^n+\int\limits_{\{v_j<u_j-\delta\}}\o_{v_j}^n]=0$ for all $\delta >0$.

Then $u_j - v_j \to 0$ in $C_X$.}

\bigskip
Recently the second named author obtained characterization of convergence in capacity for bounded $\o$-psh functions (Theorem 2.1 in [H1]). Our next result provides such a characterization in the class $\mathcal E(X,\o)$:

{\bf 2.3. Theorem.} {\sl Let $u_j\in\mathcal E(X,\o )$ and $u\in PSH (X,\o )$. Then the following three statements are equivalent:

i) $u\in\mathcal E(X,\o )$ and $u_j\to u$ in $C_X$;

ii) $\o_{u_j}^n\ll C_X$ uniformly for $j\geq 1$ and $u_j\to u$ in $C_X$;

iii) $\o_{u_j}^n\ll C_X$ uniformly for $j\geq 1$, $\varlimsup\limits_{j\to\infty} u_j\leq u$ and $\lim\limits_{j\to\infty}\int\limits_{\{u_j<u- \delta\}}\o_{u_j}^n=0$ for all $\delta >0$.}

{\sl Proof.} i) $\Leftrightarrow$ ii) This equivalence is the content of Proposition 1.11.

ii) $\Rightarrow$ iii) This implication is trivial.

iii) $\Rightarrow$ i) This is a direct application Theorem 2.1, coupled with Proposition 1.11.

\bigskip
Recall that the weak convergence of Monge-Amp\`ere measures does not imply the weak convergence of the corresponding functions (and vice versa). Thus such convergence problems are considered with additional restrictions. In [CK2] good convergence properties were obtained (in the setting of domains in $\bf C^n$) under the assumption that all the Monge-Amp\`ere measures are dominated by a fixed measure vanishing on pluripolar sets. In the K\"ahler setting the corresponding problem was studied by Xing in [Xi2]. He worked in a subclass  $\mathcal E^1(X,\o)\subset\mathcal E(X,\o)$ however all his arguments can be applied in the whole $\mathcal E(X,\o)$ provided one has an uniqueness (modulo additive constant) for the solutions of the Monge-Amp\`ere equation
$$ \o_{\va}=\mu,\ \va\in\mathcal E(X,\o),$$
where $\mu$\ is any positive Borel measure satisfying $\int\limits_X d\mu=\int\limits_X \o^n$ and vanishing on pluripolar sets. This uniqueness statement was recently obtained in [Di3]. Thus coupling Xing's arguments together with Theorem 2.2 one gets the following stability theorem:

{\bf 2.4. Theorem.} {\sl Let $u_j\in\mathcal E(X,\o )$ and $u\in\text{PSH} (X,\o )$. Assume that $\omega_{u_j}^n\leq d\mu$ for some measure $d \mu\ll C_X$ (equivalently: $\mu$\ does not charge pluripolar sets). Let also $\sup\limits_Xu_j=\sup\limits_Xu$. Then the following three statements are equivalent:

i) $u_j\to u$ in $C_X$.

ii) $u_j\to u$ in $L^1$.

iii) $\omega_{u_j}^n\to \omega_u^n$ weakly.}

The implication ii) $\Leftrightarrow$ iii) is essentially due to Xing. In particular Corollary 1 in [Xi2] (and Proposition 1.11) yields ii) $\Rightarrow$ iii). For the other implication one can proceed exactly as in Theorem 8 in [Xi2]: one can extract a subsequence from $u_j$ convergent in $L^1$\ to some function $v$. It follows that $\o_u^n=\o_v^n$, so by uniqueness (now in $\mathcal E(X,\o)$) $u=v$. Thus, since any subsequence has subsequence convergent in $L^1$\ to $u$\ one concludes that $u_j\to u$\ in $L^1$.

Since implication i) $\Rightarrow$ ii) is trivial, we only have to show that ii) implies i). But we know that $\o_{u_j}^n\ll C_X$ (since all the measures are dominated by $\mu$). By [Ce3] (see also Lemma 1.4 in [CK2]) we know that
$$\forall M\in{\bf R} \ \  \lim_{j\to\infty}\int\limits_X\max(u_j,M)-\max(u, M)d\mu=0.$$
Thus we obtain
$$\int\limits_{\{u_j<u-\delta\}}\o_{u_j}^n+\int\limits_{\{u<u_j-\delta\}}\o_{u}^n\leq\int\limits_{\{u_j<M\}\cup\{u<M\}}d\mu+\frac1{\delta}\int\limits_{X}|\max(u_j,M)-\max(u, M)|d\mu.$$
If $M$\ is sufficiently negative by the assumption $\mu\ll C_X$\ and proposition 1.7 we obtain that the first term on the right hand side can be made arbitrarily small (independently of $j$). Finally we obtain
$$\lim_{j\to\infty}[\int\limits_{\{u_j<u-\delta\}}\o_{u_j}^n+\int\limits_{\{u<u_j-\delta\}}\o_{u}^n]=0.$$
Thus by theorem 2.2 $u_j\to u$\ in $C_X$.

\bigskip
Our next result concerns a new type of stability. Instead of controlling the deviation of $u_j$\ from $v_j$ we impose an assumption on the variation of their Monge-Amp\`ere measures. Such a theorem is new even for bounded $\o$-psh functions.

{\bf 2.5. Theorem.} {\sl Let $u_j,v_j\in\mathcal E(X,\o )$ be such that

i) $\lim\limits_{j\to\infty} |\sup\limits_X u_j - \sup\limits_X v_j| =0$;

ii) $\o_{u_j}^n+\o_{v_j}^n\ll C_X$ uniformly for $j\geq 1$;

iii) $\lim\limits_{j\to\infty}\int\limits_X ||\o_{u_j}^n-\o_{v_j}^n||=0$.

Then $u_j - v_j \to 0$ in $C_X$.}

Before we start the proof we state an auxiliary lemma. To authors' knowledge such an estimate is new and might be of independent interest.

{\bf 2.6. Lemma.} {\sl Let $u,v\in\mathcal E(X,\o )$. Then}
$$\int\limits_X ||\o_u^k\w\o_v^{n-k} - \o_u^n ||\leq 2 [\int\limits_X ||\o_u^n - \o_v^n ||]^{\frac {n} {n-k}},$$

{\sl Proof.} Set 
$$\mu=\frac 1 2 [ \o_{u}^n+\o_{v}^n ].$$
We choose $f,g\in L^1 (d \mu), f,g\geq0$ such that
$$\o_{u}^n=fd\mu,\ \o_{v}^n=gd\mu.$$
By Theorem 2.1 in [Di3] and H\"older inequality we get
$$\aligned&\int\limits_E\o_{u}^k\w\o_{v}^{n-k}-\int\limits_E\o_{u}^n\geq\int\limits_Ef^{\frac k n}g^{\frac {n-k} n}d \mu-\int\limits_E fd \mu=\int\limits_E f^{\frac k n}[g^{\frac {n-k} n} - f^{\frac {n-k} n}]d \mu\\
&\geq -\int\limits_E f^{\frac k n}|g^{\frac {n-k} n} - f^{\frac {n-k} n}|d \mu\geq -\int\limits_E f^{\frac k n}|g - f|^{\frac {n-k} n}d \mu\geq -[\int\limits_E f d\mu ]^{\frac k n} [\int\limits_E |g - f| d\mu ]^{\frac {n-k} n}\\
&\geq -[\int\limits_E |g - f| d\mu ]^{\frac {n-k} n}\geq -[\int\limits_X ||\o_{u}^n - \o_{v}^n||]^{\frac {n-k} n}
\endaligned$$
for any Borel set $E\subset X$. Similarly we get
$$\int\limits_{X\backslash E}\o_{u}^k\w\o_{v}^{n-k}-\int\limits_{X\backslash E}\o_{u}^n\geq -[\int\limits_X ||\o_{u}^n - \o_{v}^n||]^{\frac {n-k} n}$$
for any Borel set $E\subset X$. Hence
$$\int\limits_{E}\o_{u}^k\w\o_{v}^{n-k}-\int\limits_{E}\o_{u}^n\leq [\int\limits_X ||\o_{u}^n - \o_{v}^n||]^{\frac {n-k} n}$$
for all Borel set $E\subset X$. Combination of these inequalities yields 
$$| \int\limits_{E}\o_{u}^k\w\o_{v}^{n-k}-\int\limits_{E}\o_{u}^n |\leq [\int\limits_X ||\o_{u}^n - \o_{v}^n||]^{\frac {n-k} n},$$
for all Borel $E\subset X$. This implies that
$$\int\limits_X ||\o_u^k\w\o_v^{n-k} - \o_u^n ||\leq 2 [\int\limits_X ||\o_u^n - \o_v^n ||]^{\frac {n} {n-k}}.$$

{\sl Proof of Theorem 2.5.}  

Case I: $u_j, v_j\in PSH (X,\omega )\cap\text{L}^\infty (X)$. Set
$$a_j=[\sup\{\int\limits_X ||\o_{u_j}^k\w\o_{v_j}^{n-k} - \o_{u_j}^m\w\o_{v_j}^{n-m}||:\ 1\leq k,m\leq n\}]^{\frac 1 4}.$$
By Lemma 2.6 we get $\lim\limits_{j\to\infty} a_j =0$. We  shall prove that 
\begin{equation}\label{2}
\lim\limits_{j\to\infty}\sup\limits_{t\in{\bf R}}\{\int\limits_{\{|u_j-v_j-t|\leq a_j\}}\o_{u_j}^n\}=1.
\end{equation}

Suppose that there exists $\e_0>0$ such that 
$$\int\limits_{\{|u_j-v_j-t|\leq a_j\}}\o_{u_j}^n\leq 1-2\epsilon_0$$
for all $t\in\bf R$, $j\geq 1$. Set
$$t_j=\sup\{t\in{\bf R}:\ \int\limits_{\{u_j<v_j+t+a_j\}}\o_{u_j}^n\leq 1-\epsilon_0\}$$
Replacing $v_j+t_j$ by $v_j$ we can assume that $t_j=0$. Then $\int\limits_{\{u_j<v_j+a_j\}}\o_{u_j}^n\leq 1-\epsilon_0$ and $\int\limits_{\{u_j\leq v_j+a_j\}}\o_{u_j}^n\geq 1-\epsilon_0$. Hence 
$$\aligned&\int\limits_{\{v_j<u_j+a_j\}}\o_{u_j}^n=1-\int\limits_{\{u_j+a_j\leq v_j\}}\o_{u_j}^n=1-\int\limits_{\{u_j\leq v_j+a_j\}}\o_{u_j}^n\\
&+\int\limits_{\{v_j-a_j<u_j\leq v_j+a_j\}}\o_{u_j}^n\leq 1-\epsilon_0.
\endaligned$$
Since $\int\limits_{\{|u_j-v_j|\leq a_j\}}\o_{u_j}^n\leq 1$ we can choose $s_j\in [-a_j+a_j^2,a_j-a_j^2]$ satisfying $$\int\limits_{\{|u_j-v_j-s_j|<a_j^2\}}\o_{u_j}^n\leq 2 {a_j}.$$
Replacing $v_j+s_j$ by $v_j$ we can assume that $s_j=0$. One easily obtains the following inequalities
$$\int\limits_{\{u_j<v_j+a_j^2\}}\o_{u_j}^n\leq 1-\epsilon_0,\ \int\limits_{\{v_j<u_j+a_j^2\}}\o_{u_j}^n\leq 1-\epsilon_0,\ \int\limits_{\{|u_j-v_j|<a_j^2\}}\o_{u_j}^n\leq 2 {a_j}.$$
By [GZ2] we can find $\rho_j\in \mathcal E(X,\omega )$, such that $\sup\limits_X\rho_j = 0$ and $\omega_{\rho_j}^n=\frac 1 {1-\epsilon_0}1_{\{u_j<v_j\}}\omega_{u_j}^n+c_j1_{\{u_j\geq v_j\}}\omega_{u_j}^n$ ($c_j\geq 0$\ is chosen so that the measure has total mass $1$). Set
$$U_j=\{(1-a_j^3) u_j < (1-a_j^3) v_j + a_j^3\rho_j\}\subset\{u_j<v_j\}.$$
By Theorem 2.1 in [Di3] we get
$$\omega_{u_j}^{n-1}\w\omega_{(1-a_j^3) v_j + a_j^3\rho_j}\geq (1-a_j^3)\omega_{u_j}^{n-1}\w\omega_{v_j}+\frac {a_j^3} {(1-\epsilon_0)^{\frac 1 n}}\omega_{u_j}^n$$
on $U_j$. By Theorem 2.3 in [Di3] we obtain
$$\aligned &(1-a_j^3)\int\limits_{U_j}\omega_{u_j}^{n-1}\w\omega_{v_j}+\frac {a_j^3} {(1-\epsilon_0)^{\frac 1 n}}\int\limits_{U_j}\omega_{u_j}^n\leq\int\limits_{U_j}\omega_{(1-a_j^3) v_j + a_j^3\rho_j}\w\omega_{u_j}^{n-1}\\ 
&\leq\int\limits_{U_j}\omega_{(1-a_j^3) u_j}\w\omega_{u_j}^{n-1}=(1-a_j^3)\int\limits_{U_j}\omega_{u_j}^n+a_j^3\int\limits_{U_j}\omega\w\omega_{u_j}^{n-1}\\
&\leq (1-a_j^3)(\int\limits_{U_j}\omega_{u_j}^{n-1}\w\omega_{v_j}+a_j^4)+a_j^3\int\limits_{U_j}\omega\w\omega_{u_j}^{n-1}.
\endaligned$$
Hence 
$$\aligned&\frac {1} {(1-\epsilon_0)^{\frac 1 n}}[\int\limits_{\{u_j\leq v_j-a_j^2\}}\omega_{u_j}^n-\int\limits_{\{\rho_j\leq -\frac {1-a_j^3}{a_j}\}}\omega_{u_j}^n]\leq \frac {1} {(1-\epsilon_0)^{\frac 1 n}}\int\limits_{U_j}\omega_{u_j}^n\\
&\leq a_j + \int\limits_{U_j}\omega\w\omega_{u_j}^{n-1}\leq a_j + \int\limits_{\{u_j<v_j\}}\omega\w\omega_{u_j}^{n-1}.
\endaligned$$
Similarly to $\rho_j$ we define $\vartheta_j\in\mathcal E(X,\o)$, such that $\sup\limits_X\vartheta_j = 0$ and $\omega_{\vartheta_j}^n=\frac 1 {1-\epsilon_0}1_{\{v_j<u_j\}}\omega_{u_j}^n+d_j1_{\{v_j\geq u_j\}}\omega_{u_j}^n$ ($d_j$\ plays the same role as $c_j$\ above).
Set
$$V_j=\{(1-a_j^3) v_j < (1-a_j^3) u_j + a_j^3\vartheta_j\}\subset\{v_j<u_j\}.$$
Similarly we get
$$\aligned&\frac {1} {(1-\epsilon_0)^{\frac 1 n}}[\int\limits_{\{v_j\leq u_j-a_j^2\}}\omega_{u_j}^n-\int\limits_{\{\vartheta_j\leq -\frac {1-a_j^3}{a_j}\}}\omega_{u_j}^n]\leq \frac {1} {(1-\epsilon_0)^{\frac 1 n}}\int\limits_{V_j}\omega_{u_j}^n\\
&\leq a_j + \int\limits_{V_j}\omega\w\omega_{u_j}^{n-1}\leq a_j + \int\limits_{\{v_j<u_j\}}\omega\w\omega_{u_j}^{n-1}.
\endaligned$$
Combination of these inequalities yields 
$$\aligned\frac {1} {(1-\epsilon_0)^{\frac 1 n}}[1-2a_j-2\int\limits_{\{\rho_j\leq -\frac {1-a_j^3}{a_j}\}}\omega_{u_j}^n]&\leq \frac {1} {(1-\epsilon_0)^{\frac 1 n}}[\int\limits_{\{|u_j-v_j|\geq a_j^2\}}\omega_{u_j}^n-2\int\limits_{\{\rho_j\leq -\frac {1-a_j^3}{a_j}\}}\omega_{u_j}^n]\\
&\leq 2 a_j + 1.
\endaligned$$
Now, by Proposition 1.7 and assumption ii), if we let $j\to\infty$  we would obtain
$$\frac {1} {(1-\epsilon_0)^{\frac 1 n}}\leq 1,$$
a contradiction.

Using (\ref{2}) we can choose $k_j\in\bf R$ such that 
$$\lim\limits_{j\to\infty}\int\limits_{\{|u_j-v_j-k_j|>a_j\}}\o_{u_j}^n=0.$$ 
Moreover from iii) we get
$$\lim\limits_{j\to\infty}\int\limits_{\{|u_j-v_j-k_j|>a_j\}}\o_{v_j}^n=0.$$ 
By Theorem 2.2 we obtain $u_j-v_j-k_j\to 0$ in $C_X$. Moreover from i) we get $\lim\limits_{j\to\infty}k_j=0$. Hence $u_j-v_j\to 0$ in $C_X$.

Case II: $u_j, v_j\in\mathcal E(X,\omega )$. We can assume, adding constants if necessary, that $\sup\limits_X u_j=\sup\limits_X v_j =0$. Choose $t_j\to +\infty$ such that $\lim\limits_{j\to\infty}\int\limits_{\{u_j\leq -t_j\}}\omega_{u_j}^n+\int\limits_{\{v_j\leq -t_j\}}\omega_{v_j}^n=0$. Using Case I for $\max (u_j,-t_j)$ and $\max (v_j,-t_j)$ we obtain $\max (u_j,-t_j)-\max (v_j,-t_j)\to 0$ in $C_X$. Therefore $u_j-v_j\to 0$ in $C_X$.

\bigskip
Our last result is a generalization of Theorem 5 from [Xi2] and Theorem 3.4 from [Ko2]. Ko\l odziej proved the statement for bounded functions, while Xing needed the additional assumption that $\forall j\geq1, v_j\geq v_0$\ for some fixed function $v_0\in\mathcal E(X,\o)$. In fact we show that this condition is superflous.

{\bf 2.7. Theorem.} {\sl Let $u_j,v_j,v\in\mathcal E(X,\o )$, $u\in\text{PSH}(X,\o )$ and $A>1$ be such that $u_j\to u$ in $L^1(X)$ and $\o_{u_j}^n\leq A \o_{v_j}^n$ for $j\geq 1$. Assume that $v_j\to v$ in $C_X$. Then $u\in\mathcal E(X,\o)$ and $u_j\to u$ in $C_X$.}

{\sl Proof.} For each $t>0$ we set
$$u_{jt}:=\max (u_j, -t),\ u_{t}:=\max (u, -t),\ v_{jt}:=\max (v_j, -t),\ v_{t}:=\max (v, -t),$$
$$T_{jt}:=\sum\limits_{k=0}^{n-1}\o_{v_{jt}}^k\w\o_{v_{t}}^{n-1-k}.$$ 
We have
$$\int\limits_{\{u_j<u- \delta\}}\o_{u_j}^n\leq A\int\limits_{\{u_j<u- \delta\}}\o_{v_j}^n.$$
The latter integral can be estimated by
\begin{align*}&\int\limits_{\{u_j<u- \delta\}}\o_{v_j}^n\leq\int\limits_{\{u_j\leq -t\}\cup\{u\leq -t\}\cup\{v_j\leq -t\}}\o_{v_j}^n+\int\limits_{\{u_{jt}<u_t- \delta\}}\o_{v_{jt}}^n\leq\int\limits_{\{u_j\leq -t\}\cup\{u\leq -t\}\cup\{v_j\leq -t\}}\o_{v_j}^{n}\\
&+\frac 1 {\delta}\int\limits_{X}|u_{t}-u_{jt}|\o_{v_{jt}}^n\leq\int\limits_{\{u_j\leq -t\}\cup\{u\leq -t\}\cup\{v_j\leq -t\}}\o_{v_j}^n+\frac 1 {\delta}\int\limits_{X}|u_{t}+\e -u_{jt}|\o_{v_{jt}}^n\\
&+\frac {\e}{\delta}\leq\int\limits_{\{u_j\leq -t\}\cup\{u\leq -t\}\cup\{v_j\leq -t\}}\o_{v_j}^n+\frac 1 {\delta}\int\limits_{X}(u_{t}-u_{jt})\o_{v_{jt}}^n+\frac {t+\e } {\delta}\int\limits_{\{u_{jt}>u_t+\e \}}\o_{v_{jt}}^n+\frac {2\e }{\delta}\\
&\leq\int\limits_{\{u_j\leq -t\}\cup\{u\leq -t\}\cup\{v_j\leq -t\}}\o_{v_j}^n+\frac 1 {\delta}\int\limits_{X}(u_{t}-u_{jt})(\o_{v_{jt}}^n-\o_{v_t}^n)\ +\frac 1 {\delta}\int\limits_{X}(u_{t}-u_{jt})\o_{v_{t}}^n
\end{align*}
\begin{align*}
&+\frac {t+\e } {\delta}\int\limits_{\{u_{jt}>u_t+\e \}}\o_{v_{jt}}^n+\frac {2\e }{\delta}=\int\limits_{\{u_j\leq -t\}\cup\{u\leq -t\}\cup\{v_j\leq -t\}}\o_{v_j}^n+\frac 1 {\delta}\int\limits_{X}(v_{t}-v_{jt})(\o_{u_{jt}}-\o_{u_t})\w T_{jt}\\
&+\frac 1 {\delta}\int\limits_{X}(u_{t}-u_{jt})\o_{v_{t}}^n+\frac {t+\e } {\delta}\int\limits_{\{u_{jt}>u_t+\e \}}\o_{v_{jt}}^n+\frac {2\e }{\delta}\leq\int\limits_{\{u_j\leq -t\}\cup\{u\leq -t\}\cup\{v_j\leq -t\}}\o_{v_j}^n\\
&+\frac 1 {\delta}\int\limits_{X}|v_{t}-v_{jt}|(\o_{u_{jt}}+\o_{u_t})\w T_{jt}+\frac 1 {\delta}\int\limits_{X}(u_{t}-u_{jt})\o_{v_{t}}^n+\frac {t+\e } {\delta}\int\limits_{\{u_{jt}>u_t+\e \}}\o_{v_{jt}}^n+\frac {2\e }{\delta}
\end{align*}
for all $t>0$ and all $\e >0$. 

Observe that the second term goes to $0$\ aj $j\to\infty$, since the intergated function converges to $0$\ in capacity and all the measures are uniformly absolutely continuous with respect to $C_X$ for $j\geq1$ (and fixed $t$).
By [Ce3] we also get 
$$\lim\limits_{j\to\infty} \int\limits_{X}(u_{t}-u_{jt})\o_{v_{t}}^n=0.$$
Note that Hartogs' lemma (see [H\"o]) yields $C_X (\{u_{jt}>u_t+\e \})\to 0$ as $j\to\infty$. Moreover, since $\omega_{v_{jt}}^n\ll C_X$ uniformly for $j\geq 1$ we obtain 
$$\lim\limits_{j\to\infty}\int\limits_{\{u_{jt}>u_t+\e \}}\o_{v_{jt}}^n=0.$$

Coupling all these and letting $j\to\infty$  we get 
$$\varlimsup\limits_{j\to\infty}\int\limits_{\{u_j<u- \delta\}}\o_{u_j}^n\leq\sup\limits_{j\geq 1}\int\limits_{\{u_j\leq -t\}\cup\{u\leq -t\}\cup\{v_j\leq -t\}}\o_{v_j}^n+\frac {2\e }{\delta},$$
for all $t>0$, $\e >0$. Letting $\e\to 0$ and $t\to +\infty$ we obtain
$$\lim\limits_{j\to\infty}\int\limits_{\{u_j<u- \delta\}}\o_{u_j}^n=0.$$
By Theorem 2.3 we obtain that $u_j\to u$ in $C_X$.

\bigskip
\noindent
\centerline{\bf REFERENCES}

\bigskip
\noindent
[Bl] Z. B\l ocki, {\sl Uniqueness and stability for the complex Monge-Amp\`{e}re equation on compact K\"ahler manifolds}, Indiana Univ. Math. J. {\bf 52} (2003), no. 6, 1697-1701.

\noindent
[BT1] E. Bedford and B. A. Taylor, {\sl The Dirichlet problem for the complex Monge-Amp\`{e}re operator}. Invent. Math. {\bf 37} (1976), 1-44.

\noindent
[BT2] E. Bedford and B. A. Taylor, {\sl A new capacity for plurisubharmonic functions}. Acta Math. {\bf 149} (1982), 1-40.

\noindent
[Ce1] U. Cegrell, {\sl Pluricomplex energy}, Acta Math. {\bf 180} (1998), 187-217.

\noindent
[Ce2] U. Cegrell, {\sl The general definition of the complex Monge-Amp\`{e}re operator}, Ann. Inst. Fourier (Grenoble) {\bf 54} (2004), 159-179.

\noindent
[Ce3] U. Cegrell, {\sl Convergence in capacity}, Technical report, Issac Newton Institute for Mathematical Sciences, 2001.

\noindent
[CK1] U. Cegrell, S. Ko\l odziej, {\sl The Dirichlet problem for the complex Monge-Amp\`ere operator: Perron classes and rotation invariant measures}, Mich. Math. J. {\bf 41} (1994), 563-569.

\noindent
[CK2] U. Cegrell, S. Ko\l odziej, {\sl The equation of complex Monge-Amp\`ere type and stability of solutions}, Math. Ann. {\bf 334} (2006), 713-729.

\noindent
[CGZ] D. Coman, V. Guedj and A. Zeriahi, {\sl Domains of definition of Monge-Amp\`{e}re operators on compact K\"ahler manifolds}, Math. Zeit. {\bf 259} (2008), 393-418.

\noindent
[Di1] S. Dinew, {\sl Cegrell classes on compact K\"ahler manifolds}, Ann. Polon. Math. {\bf 91} (2007), 179-195.

\noindent
[Di2] S. Dinew, {\sl An inequality for mixed Monge-Amp\`{e}re measures}, Math. Zeit. {\bf 262} (2009), 1-15.

\noindent
[Di3] S. Dinew, {\sl Uniqueness in $\mathcal E(X,\o )$}, J. Funct. Anal. {\bf 256} (2009), 2113-2122.

\noindent
[GZ1] V. Guedj and A. Zeriahi, {\sl Intrinsic capacities on compact K\"ahler manifolds}, J. Geom. Anal. {\bf 15} (2005), no. 4, 607-639.

\noindent
[GZ2] V. Guedj and A. Zeriahi, {\sl The weighted Monge-Amp\`{e}re energy of quasiplurisubharmonic functions}, J. Funct. Anal. {\bf 250} (2007), 442-482.

\noindent
[H1] P. Hiep, {\sl On the convergence in capacity on compact K\"ahler manifolds and its applications}, Proc. Amer. Math. Soc. {\bf 136} (2008), 2007-2018.

\noindent
[H2] P. Hiep, {\sl Convergence in capacity}, Ann. Polon. Math. {\bf 93} (2008), 91-99.

\noindent
[H\"o] L. H\"ormander, {\sl Notions of Convexity}, Progess in Mathematics {\bf 127}, Birkh\"auser, Boston, (1994).

\noindent
[KH] N. Khue and P. Hiep, {\sl A comparison principle for the complex Monge-Amp\`{e}re operator in Cegrell's classes and applications}, to appear in Trans. Amer. Math. Soc.

\noindent
[Ko1] S. Ko\l odziej, {\sl The Monge-Amp\`{e}re equation on compact K\"ahler manifolds}, Indiana Univ. Math. J. {\bf 52} (2003), 667-686.

\noindent
[Ko2] S.Ko\l odziej, {\sl The set of measures given by bounded solutions of the complex Monge-Amp\`ere equation on compact K\"ahler manifolds}, J. London Math. Soc. {\bf 72} (2) (2005), 225-238.

\noindent
[Xi1] Y. Xing, {\sl Continuity of the complex Monge-Amp\`{e}re operator}, Proc. Amer. Math. Soc. {\bf 124} (1996), 457-467.

\noindent
[Xi2] Y. Xing {\sl Continuity of the Complex Monge-Amp\`{e}re Operator on Compact K\"ahler Manifolds}, to appear in Math. Zeit. (http://www.arxiv.org/).

\bigskip
Jagiellonian University,

Mathematics Institute,

30-348 \L ojasiewicza 6,

Krak\'ow, Poland

\tt slawomir.dinew@im.uj.edu.pl

\rm

\bigskip
Department of Mathematics

University of Education (Dai hoc Su Pham  Ha Noi)

CauGiay, Hanoi, Vietnam

\tt  phhiep\underline{\text{ }}vn@yahoo.com
\rm
\end{document}